\documentclass[letterpaper,10pt]{article}

\usepackage[english]{babel}
\usepackage{amssymb}
\usepackage{latexsym}
\usepackage{amsfonts}
\usepackage{amsmath}
\usepackage{amscd}
\usepackage[dvips]{graphicx}

\newcommand{\FAT}[1]{\mbox{{$\mathbb{#1}$}}}

\newcommand{\HSM}[1]{\mathcal{#1}}
\newcommand{\EE}{\FAT{E}}

\newcommand{\NN}{\FAT{N}}

\newcommand{\ZZ}{\FAT{Z}}

\newcommand{\RR}{\FAT{R}}

\newcommand{\HH}{\HSM{H}}

\newcommand{\THEN}{\,\Rightarrow\,}
\newcommand{\IFF}{\,\Leftrightarrow\,}

\newcommand{\ep}{\mbox{$\quad\blacksquare$}}

\newcommand{\proof}[1]{\textbf{Proof #1: }}
\newcommand{\minus}{\smallsetminus}
\newcommand{\cl}[1]{\overline{#1}}
\renewcommand{\int}[1]{\stackrel{\circ}{#1}}
\newcommand{\bd}{\partial_\infty}
\renewcommand{\ae}{\stackrel{ae}{=}}

\newcommand{\med}{\mathbf{med}}

\newcommand{\sh}[1]{\mathrm{sh}\left(#1\right)}

\newcommand{\shcirc}[1]{\mathrm{sh}^\circ\left(#1\right)}
\newcommand{\gsh}[1]{\mathrm{gsh}\left(#1\right)}

\newcommand{\xhat}{X\cup\bd X}
\newcommand{\hgt}[1]{\left|#1\right|}

\newtheorem{thm}{Theorem}[section]
\newtheorem{prop}[thm]{Proposition}
\newtheorem{lemma}[thm]{Lemma}
\newtheorem{defn}[thm]{Definition}

\newtheorem{cor}[thm]{Corollary}
\newtheorem{remark}[thm]{Remark}

\begin{document}
\title{Local finiteness of cubulations and CAT(0) groups}
\author{Dan P. Guralnik\\ Vanderbilt University}
\maketitle
\begin{abstract} Let $X$ be a geodesic space and $G$ a group acting geometrically on $X$. A discrete halfspace system of $X$ is a set $\HH$ of open halfspaces closed under $h\mapsto X\minus\cl{h}$ and such that every $x\in X$ has a neighbourhood intersecting only finitely many walls of $\HH$. Given such a system $\HH$, one uses the Sageev-Roller construction to form a cubing $C(\HH)$. When $\HH$ is invariant under $G$ we have:

\vspace{1.5mm}
\noindent\textbf{Theorem A. }{\it $X$ has a $G$-equivariant quasi-isometric embedding into $C(\HH)$.}

\vspace{1.5mm}
The basic questions about $C(\HH)$ are: when are all cubes in $C(\HH)$ finite-dimensional? when is $C(\HH)$ finite dimensional? when is it proper? when is $C(\HH)$ $G$-co-compact (and hence $G$ is biautomatic, by a result of Niblo and Reeves)?

These questions were answered by Niblo-Reeves, Williams and Caprace for the case of Coxeter groups $(W,R)$ acting on their Davis-Moussong complexes, with elements of $\HH$ being the halfspaces defined by reflections. A significant role was played by the `parallel walls property' of Coxeter groups, conjectured by Davis and Shapiro and proved by Brink and Howlett. It thus becomes natural to ask these questions whenever $X$ is a CAT(0) space carrying a geometric action by a group $G$.

In this paper we show that, when $\HH$ has bounded chambers, the parallel walls property is equivalent to a condition we call \emph{uniformness}, regarding the quality of approximation of boundary points by walls of $\HH$. Uniformness, as opposed to the parallel walls property, involves no explicit bounds. We prove:

\vspace{1.5mm}
\noindent\textbf{Theorem B. }{\it Let $G$ be a group acting geometrically on a geodesic space $X$ and suppose $\HH$ is a discrete $G$-invariant halfspace system in $X$. If $\HH$ is uniform, then $C(\HH)$ is proper (locally-finite). In particular, $C(\HH)$ does not contain infinite-dimensional cubes.}
\end{abstract}
\section{Introduction}
Discrete systems of walls are a very old object. They arise naturally in geometry in connection with discontinuous isometric actions of groups on metric spaces. For example, in the model geometries of constant curvature, one uses halfspaces defined using the relevant distance function in order to construct fundamental domains for any given action (Dirichlet domains).

The study of these systems in their own right has gained momentum after the discovery by Sageev, followed by Roller, of a duality between such systems and non-positively curved cube complexes. When viewed as a partially-ordered set (ordered by inclusion) with a complementation operator (switching every halfspace with its complementary counterpart), such a system $\HH$ in a space $X$ gives rise to a cube complex $C(\HH)$ whose vertices correspond to the principal ultrafilters on $\HH$. Even the most basic natural examples -- those of triangle Coxeter groups in the Euclidean plane -- show that it is very hard to control the dimension and the growth properties of $C(\HH)$.

Except for the case when $X$ is itself a non-positively curved (piecewise-Euclidean) cube complex and $\HH$ is its natural system of halfspaces, the best studied situation in this respect is that of Coxeter groups of finite rank. To any Coxeter system $(W,R)$ corresponds its Davis-Moussong complex $X=M(W,R)$, carrying a natural piecewise-Euclidean CAT(0) metric, on which the reflections of $(W,R)$ act as actual reflections (orientation-reversing isometries of order $2$, with nowhere-dense, convex fixed-point sets, separating the space into two convex components). Thus, the fixed point sets of reflections may serve as a system of walls in $X$, giving rise to a halfspace system $\HH$.  Combining results of Brink and Howlett \cite{[BriHow]} and Niblo-Reeves \cite{[NibRee1]}, one sees the following pattern of ideas:
\begin{enumerate}
	\item There is a bound on the dimension of $C(\HH)$, and this bound is used for the proof of local finiteness (!).
	\item There is a bound on the distance of any point $p\in X$ from any wall of $\HH$ not separated from $p$ by other walls -- this is the `parallel walls property';
	\item The above bounds are related through properties of the root system corresponding to $(W,R)$;
	\item Combining these bounds using cancellation properties of $(W,R)$ one obtains the local finiteness of $C(\HH)$.
\end{enumerate}
A later work by Williams \cite{[Wil]} addressed the co-compactness problem for Coxeter groups: when does $G$ act co-compactly on $C(\HH)$? Williams managed to provide a partial answer, and the discussion was finished by Caprace \cite{[Caprace]}, who proved Williams' conjecture that the action of $W$ on $C(\HH)$ is co-compact iff $W$ does not contain an Euclidean triangle subgroup. Caprace also obtained a uniform bound on the degrees of vertices of $C(\HH)$, including the non co-compact case, strengthening the `parallel walls property' for Coxeter groups: using the Caprace bounds one is able to tell how far should a point of $M(W,R)$ lie from a given wall so that there is a prescribed number of intermediate walls.

The current work arose as part of an effort to understand the extent to which similar results remain true for general groups acting on CAT(0) spaces, as halfspace systems -- as well as their `cousins', spaces with walls -- provide one of the main tools for constructing group actions on cubings, and one wants to have the latter as tame as possible.

Our main result generalizes and strengthens the local-finiteness result of Niblo and Reeves in the following manner: given a group $G$ acting on a CAT(0) space $X$ and preserving a discrete halfspace system $X$, we are able to prove $\HSM{C}=C(\HH)$ is locally finite, when our starting point is that 
\begin{enumerate}
	\item we replace the `parallel walls property' by a seemingly weaker assumption about how boundary points are approximated by walls of $\HH$,
	\item we assume $\HH$ induces bounded chambers on $X$,
	\item we assume $G$ acts co-compactly on $X$.
\end{enumerate}
Halfspace systems satisfying 1.-3. are said to be {\it uniform}, by analogy with requirements of conical convergence arising in Kleinian groups and, more generally in the theory of relatively-hyperbolic groups. The second condition is, roughly speaking, that every ideal boundary point $\xi\in\bd X$ is `well-approximated' by walls of $\HH$ cutting the representative rays of the class $\xi$ transversely. For any such $\xi$, its cone-neighbourhoods are then exhausted by compact sets arising as their intersections with descending sequences of halfspaces from $\HH$. Our main result is:

\vspace{2mm}
\noindent\textbf{Theorem B. }{\it Let $G$ be a group acting geometrically on a CAT(0) space $X$ and suppose $\HH$ is a halfspace system in $X$ invariant under $G$. If $\HH$ is uniform, then $C(\HH)$ is proper.}

\vspace{2mm}
Essential ingredients of the proof are the co-compactness {\it and} the properness of the action, as well as the resulting compactness of the boundary $\bd X$.\\

As by-products of our technique we prove that, for a halfspace system $\HH$, uniformness is, in fact, equivalent to having bounded chambers together with the parallel walls property. We also show that when $C(\HH)$ is $G$-co-compact, $\HH$ must satisfy the strong parallel walls property.\\

Another related result is 

\vspace{2mm}
\noindent\textbf{Theorem A. }{\it Let $G$ be a group acting geometrically on a geodesic space $X$ and suppose $\HH$ is a discrete $G$-invariant halfspace system in $X$. Then $X$ has a $G$-equivariant quasi-isometric embedding into $C(\HH)$.}

\vspace{2mm}
A similar result was claimed in \cite{[NibRee1]}, but only for Coxeter groups acting on their Davis-Moussong complex, and without proof. We provide a simple proof of this fact, which is later used in this paper to derive the more subtle metric properties of uniform systems. For example, we sharpen the parallel walls property by showing that the number of walls separating a point $x$ from a halfspace $h\in\HH$ must grow linearly with $d(x,h)$.\\

This paper is organized as follows. Section \ref{section:prelim} discusses some preliminary notions: visual boundaries of CAT(0) spaces, halfspaces and cubings; the Sageev-Roller duality is discussed in some detail to lay the technical groundwork to our method. Section \ref{section:halfspaces} briefly discusses the geometry of halfspace systems in geodesic and CAT(0) spaces. Section \ref{section:cubing dual to a halfspace system} discusses the geometry of the `embedding' of a space $X$ in the cubing $C(\HH)$ dual to a halfspace system $\HH$ in $X$. Finally, in section \ref{section:geometry of uniform systems} we prove the main results and discusses the geometry of uniform halfspace systems and parallel walls properties.

\paragraph{Acknowledgements.} The author is grateful to Romain Tessera and Mike Mihalik for their valuable advice regarding improvements to the exposition.

\section{Preliminaries}\label{section:prelim}
\subsection{Some CAT(0) geometry.}\label{prelim:geometry}
From now on let $(X,d)$ be a fixed proper CAT(0) space. Before we investigate halfspaces, let us recall some notions from CAT(0) geometry.
\paragraph{Visual CAT(0) boundaries. } A good reference for the content of this paragraph is chapter II of \cite{[BH]}. The space $X$ has a natural compactification by its visual boundary: we let $\bd X$ denote the set of asymptoticity classes of geodesic rays in $X$. Recall that two geodesic rays $\gamma,\gamma':[0,\infty)\to X$ are said to be asymptotic, if they fellow-travel. In CAT(0) geometry this is equivalent to their images lying at finite Hausdorff distance from each other. Another fact allowing to construct a compact topology on the space $\xhat$ is the existence of projections in $X$: for every closed convex subspace $F$ of $X$ there is a canonical map $pr_F:X\to F$ mapping any $x\in X$ to the unique point on $F$ lying at a minimal distance to $x$; one then uses this fact to construct \emph{cone neighbourhoods} as follows --
\begin{displaymath}
	U_{x_0,\xi}(R,\epsilon)=\left\{x\in\xhat\,\left|
		d\left(pr_{B(x_0,R)}(x),\gamma(R)\right)<\epsilon
	\right.\right\}\,,
\end{displaymath}
where $x_0\in X$, $\xi\in\bd X$, $R,\epsilon>0$, and we note that the projection of $X$ onto any closed ball $B(x_0,R)$ extends naturally to $\xhat$.

Fixing $x_0$, the topology on $\xhat$ generated by the metric $d$ on $X$ and the set of all $U_{x_0,\xi}(R,\epsilon)$ is called \emph{the cone topology} and is known to be independent of the choice of basepoint. $\bd X$ is called the \emph{visual boundary of $X$} when endowed with this topology. $\bd X$ is compact in the cone topology whenever $X$ is proper.

It is sometimes beneficial to consider the set of accumulation points of a subset $A$ of $X$ in $\xhat$:
\begin{defn}[ideal boundary of a subspace] the ideal boundary $\bd A$ of a subspace $A$ of $X$ equals the intersection of $\bd X$ with the closure of $A$ in $\xhat$ relative to the cone topology.
\end{defn}

\paragraph{Halfspaces and walls. } Since $X$ is uniquely-geodesic, it makes perfect sense to consider decompositions of $X$ into pairs of complementary halfspaces:
\begin{defn} A \emph{halfspace} $h\subset X$ in $X$ is a non-empty open convex subset such that $h^\ast:= X\minus\cl{h}$ is also convex. The intersection $\cl{h}\cap\cl{h^\ast}$ will be called \emph{the wall associated with $h$}, and denoted by $W(h)$; if $S$ is a set of halfspaces, then $W(S)$ will denote the set of walls $W(h)$ for $h\in S$. The sets $\varnothing, X$ are, by definition, the \emph{trivial} halfspaces of $X$. 
\end{defn}
Here we list some properties of halfspaces in CAT(0) spaces, emphasizing relations with the boundary. Let $h$ be a fixed halfspace. The most important observation about $h$ is that $\cl{h}$ is then a complete CAT(0) space with respect to the metric induced from $X$. Then so is the corresponding wall $W(h)=\cl{h}\cap\cl{h^\ast}$. Thus, it makes sense to consider the visual boundaries of $\cl{h}$ and $W(h)$, which may be constructed by computing their respective closures in $\xhat$, and intersecting those with $\bd X$. An easy consequence of these observations is:
\begin{lemma} For any halfspace $h$ in $X$ we have $\bd W(h)=\bd h\cap\bd h^\ast$.
\end{lemma}
Another consequence of the convexity of halfspaces is:
\begin{lemma} For any halfspace $h$ in $X$ and $\xi\in\bd X$, if $\xi$ is an interior point of $\bd h$, then, for any ray $\gamma\in\xi$ we have $d(\gamma(t),h^\ast)\to\infty$ as $t\to\infty$. 
\end{lemma}
These two facts will be used without reference in what follows.

\subsection{Cubings.}\label{prelim:cubings}
Let $Q^d$ denote the unit Euclidean $d$-cube, and let $\partial Q^d$ denote its $(d-1)$-skeleton. In \cite{[Sa]}, Sageev shows how a cubing $C$ is reconstructed from the metric structure on its $1$-skeleton $C^1$ (with respect to the combinatorial metric): starting with the standard geometric realization of $C^1$, one glues copies of $Q^2$ (using isometries) to fill-in each $4$-cycle in $C^1$, resulting in a square $2$-complex $C^2$; then one proceeds inductively by gluing a copies of $Q^d$ onto $C^{d-1}$ to fill-in every copy of $\partial Q^d$.

Roller in \cite{[Rol]} provides a characterization of the graphs arising as $1$-skeleta of cubings, showing that all such graphs arise as duals to certain ordered structures that he calls {\it poc-sets} (i.e., posets with complementation). Sageev's construction of cubings for multi-ended group pairs (also in \cite{[Sa]}) is a special case of this general principle. Roller also shows that what facilitates Sageev's construction is the fact that the $1$-skeleton of a cubing has the structure of a discrete median algebra.

\subsubsection{discrete poc-sets and their duals.}\label{subsubsection:ultrafilters} Since we will be using Roller's characterization, we provide all the necessary terminology here.
\begin{defn}[poc-set, nesting, transversality]\label{defn:poc-set terminology} A poc-set $(H,\leq,\ast)$ is a partially-ordered set $(H,\leq)$ with a minimum element $0$ and an order-reversing involution $h\mapsto h^\ast$ satisfying the requirement that for all $h\in
H$, if $h\leq h^\ast$ then $h=0$. 
\begin{itemize}
	\item[-] the elements $0,0^\ast$ are the \emph{trivial} elements of $H$, while all other elements of $H$ are \emph{proper}.
	\item[-] the poc-set $(H,\leq,\ast)$ is said to be \emph{discrete}, if, for every pair of proper elements $a,b\in H$, the order interval $[a,b]=\{h\in H\,|\,a\leq h\leq b\}$ is finite.
	\item[-] two elements $h,k\in H$ are said to be nested (resp. transverse), -- denoted here with $h\| k$ (resp. $h\pitchfork k$) -- if one (resp. none) of the relations $h\leq k,\,h^\ast\leq k,\,h\leq k^\ast,\,h^\ast\leq k^\ast$ holds. A subset $S\subseteq H$ is nested (resp. transverse) if all its elements are pairwise nested (resp. transverse).
	\item[-] a poc-set $(H,\leq,\ast)$ is said to be of dimension $\omega$, if it contains no infinite transverse subset.
\end{itemize}
\end{defn}
Given a discrete poc-set $H$, we consider the Stone space $2^H$ (endowed with the Tychonoff topology). On $2^H$ one has the following {\it median operation}:
\begin{displaymath}
	\med(\alpha,\beta\,\gamma)=(\alpha\cap\beta)\cup(\beta\cap\gamma)\cup(\gamma\cap\alpha)
\end{displaymath}

One restricts the median operation to the set of {\it ultrafilters on $H$}:
\begin{defn}\label{defn:ultrafilter} Suppose $(H,\leq,\ast)$ is a discrete poc-set. An ultrafilter $\alpha$ on $H$ is a subset of $H$ satisfying:
\begin{description}
    \item[$\mathbf{(UF1)}$] for all $h\in H$, either $h\in\alpha$ or $h^\ast\in\alpha$, but not both;
    \item[$\mathbf{(UF2)}$] for all $h,k\in\alpha$, the relation $h\leq k^\ast$ is prohibited.
\end{description}
The space of all ultrafilters on $H$ will be denoted by $H^\circ$.

A collection $\alpha\subset H$ satisfying $\mathbf{(UF2)}$ is called a \emph{filter base}.
\end{defn}
\textbf{Remark: } One may employ Zorn's lemma to show that any filter base is contained in an ultrafilter (see \cite{[Rol]}, 3.4(iii)).\\

The space $H^\circ$ turns out to be a closed subspace of $2^H$, which is also a median subalgebra. $H^\circ$ admits a natural distance function (we allow the value $\infty$):
\begin{displaymath}
	\Delta(\alpha,\beta)=\frac{1}{2}\left|\alpha\vartriangle\beta\right|\,.
\end{displaymath}
Note that $(\alpha\minus\beta)^\ast=\beta\minus\alpha$, by $\mathbf(UF1)$, so that 
\begin{displaymath}
	\Delta(\alpha,\beta)=\left|\alpha\minus\beta\right|=\left|\beta\minus\alpha\right|\,.
\end{displaymath}
$\Delta$ induces a natural equivalence relation $(\ae)$ on $H^\circ$: we say that $\alpha,\beta\in H^\circ$ are {\it almost-equal} (denoted $\alpha\ae\beta$), if $\Delta(\alpha,\beta)$ is finite. Every almost-equality class $\Sigma$ is turned into a graph $\Gamma_\Sigma$ with $V\Gamma_\Sigma=\Sigma$ by joining $\alpha,\beta\in\Sigma$ by an edge if and only if $\Delta(\alpha,\beta)=1$. It is an important fact that the restriction of $\Delta$ to $\Sigma$ coincides with the combinatorial metric on the (connected) graph $\Gamma_\Sigma$.

\subsubsection{Almost-equality classes as graphs.} For each almost-equality class $\Sigma$ of $H^\circ$, the restriction of $\Delta$ to $\Sigma$ is an actual metric, and one defines {\it intervals} by: 
\begin{displaymath}
	[\alpha,\beta]=\left\{\mu\in\Sigma\,\left|
		\Delta(\alpha,\beta)=\Delta(\alpha,\mu)+\Delta(\mu,\beta)
	\right.\right\}
\end{displaymath}
It turns out that this notion of intervals (for $\alpha\ae\beta$) coincides with the more general notion of intervals defined by the median algebra structure:
\begin{displaymath}
	[\alpha,\beta]=\left\{\med(\alpha,\beta,\mu)\,\left|\mu\in H^\circ\right.\right\}\,,
\end{displaymath}
and, in fact, one has the following equality:
\begin{displaymath}
	\{\med(\alpha,\beta,\gamma)\}=[\alpha,\beta]\cap[\beta,\gamma]\cap[\gamma,\alpha].
\end{displaymath}
When restricted to an almost-equality class, this equality produces a geometric interpretation for the median operation: for every triple of vertices, there exists a unique vertex lying on the intersection of the sides of {\it a} geodesic triangle with the given vertices.\\

It will be important for us to describe the local structure of an almost-equality class as a graph more precisely. For any $\alpha,\beta\in H^\circ$, if $\Delta(\alpha,\beta)=1$ then one can write $\alpha\minus\beta=\{a\}$ for some $a\in H$, which means that $\beta=(\alpha\minus a)\cup\{a^\ast\}$. We denote
\begin{equation}\label{eqn:moving one step in H circ}
	[\alpha]_a=(\alpha\minus a)\cup\{a^\ast\}\,.
\end{equation}
It is clear that, for $\alpha\in H^\circ$ one has $[\alpha]_a\in H^\circ$ if and only if $a\in\min(\alpha)$ when $\alpha$ is viewed as a subset of $H$ with the induced partial ordering. Thus, $\min(\alpha)$ parametrizes the vertices adjacent to $\alpha$ in $H^\circ$.\\

Sageev has made the observation that, if $A=\{a_1,\ldots,a_d\}$ is a transverse subset of $\min(\alpha)$ then for any permutation $\sigma\in S_n$ the ultrafilter $\left[[\alpha]_{a_{\sigma(1)}\ldots}\right]_{a_{\sigma(d)}}$ is well-defined and independent of $\sigma$. It follows that one can define $[\alpha]_B$ for all $A\subseteq B$ and the set of all such $[\alpha]_B$ forms the $1$-skeleton of a $d$-dimensional cube of the cubing $C(\Sigma)$ associated with the graph $\Gamma_\Sigma$, where $\Sigma$ is the almost-equality class of $\alpha$. Clearly, every $d$-dimensional cube in any of the cubings associated with almost-equality classes of $H^\circ$ arises in this way.

\subsubsection{Properness and dimension.}\label{subsubsection:proper cubings}
For details regarding general median algebras -- the reader is referred to \cite{[Rol]}. Two important classes of ultrafilters arise in the study of the duality between discrete poc-sets and discrete median algebras.
\begin{defn} Let $(H,\leq,\ast)$ be a discrete poc-set. An ultrafilter $\alpha\in H^\circ$ is said to be {\it well-founded}, if for every $a\in\alpha$, the set of $b\in\alpha$ satisfying $b\leq a$ is finite. More generally, $\alpha\in H^\circ$ is said to be {\it principal}, if it contains no infinite strictly-descending chain.
\end{defn}
\begin{remark}
Clearly, a well-founded ultrafilter is principal. It is also obvious from the definition that every element $a$ of a principal ultrafilter $\alpha$ has some $a_0\in\min(\alpha)$ satisfying $a_0<a$.
\end{remark}
We need the following results regarding well-founded ultrafilters:
\begin{prop}[\cite{[Rol]}, proposition 9.3] A poc set $P$ has a well-founded ultrafilter $\sigma$ if and only if there exists a discrete median algebra $M$ such that $P\cong M^\circ$. In fact, $M$ may be chosen to be the almost-equality class of $\sigma$.
\end{prop}
\begin{prop}[\cite{[Rol]}, proposition 5.9] The image of the double dual map $M\to M^{\circ\circ}$ is dense in $M^{\circ\circ}$.
\end{prop}
Combining this with the fact that a discrete median algebra is the vertex-set of a median graph (which follows from \cite{[Rol]}, proposition 2.16), it follows that if a discrete poc-set $H$ has a well-founded ultrafilter $\alpha$, then the almost-equality class of $\alpha$ is Tychonoff-dense in $H^\circ$.\\

For every ultrafilter $\alpha$ and every proper $h\in H$ let $n(\alpha,h)$ denote the number of $x\in\alpha$ satisfying $x\leq h$. If now $a\in\min(\alpha)$ then for every $a^\ast\neq h\in[\alpha]_a$ we have $h\in\alpha$ and $n([\alpha]_a,h)\leq n(\alpha,h)+1$ while for $h=a^\ast$ we have $n([\alpha]_a,h)=1$. In particular, if $\alpha$ is well-founded then every element of the almost-equality class of $\alpha$ is well-founded. A similar argument shows that principal ultrafilters also `come in almost-equality classes'.\\

Suppose now that $H$ has a well-founded ultrafilter $\sigma$ whose degree (as a vertex in the graph $\Gamma=\bigsqcup\Gamma_\Sigma$, $\Sigma$ ranging over the almost-equality classes of $H^\circ$) is finite. Equivalently $\min(\sigma)$ is finite. Let $\Sigma_0$ denote the almost-equality class of $\sigma$.

Suppose $\pi\in H^\circ\minus\Sigma_0$, so that $\pi\minus\sigma$ is infinite. Consider two situations:
\begin{enumerate}
	\item If $\pi$ is well-founded, then every
\begin{displaymath}
	h\in\pi\minus\sigma=\pi\cap\sigma^\ast
\end{displaymath}
lies in an interval of the form $[a,b^\ast]$ for $a\in\min(\pi)$ and $b\in\min(\sigma)$. Since $\min(\sigma)$ is finite, we conclude there exists a $b_0\in\min(\sigma)$ such that the set of elements $h\in\pi\minus\sigma$ satisfying $h<b_0^\ast$ is infinite. However, $b_0^\ast$ must be in $\pi$, and we obtain a contradiction to the well-foundedness of $\pi$. Thus, every well-founded ultrafilter on $H$ is almost-equal to $\sigma$.
	\item More generally, if $\pi$ is principal, then, once again, an infinite subset of $\pi\minus\sigma$ is covered by intervals of the form $[a,b_0^\ast]$ with $a\in\min(\pi)$ and $b_0$ a fixed element of $\min(\sigma)$. Suppose $a_1,a_2\in\min(\pi)$ satisfy $a_1,a_2<b_0^\ast$ and $a_1^\ast\leq a_2$; then $b_0<a_1^\ast\leq a_2<b_0^\ast$, contradicting $b_0$ being proper; since $a_1,a_2$ are incomparable and their both lying in $\pi$ makes $a_1\leq a_2^\ast$ impossible, we conclude that $a_1$ and $a_2$ are transverse. Thus, the assumption that $\pi\not\ae\sigma$ implies $\min(\pi)$ contains an infinite transverse set.
\end{enumerate}
We have proved:
\begin{prop}\label{prop:principal class exists} Let $(H,\leq,\ast)$ be a discrete poc-set. If $H$ has a well-founded ultrafilter $\sigma$ with $\min(\sigma)$ finite, then:
\begin{enumerate}
	\item The well-founded ultrafilters of $H$ constitute a dense almost-equality class $\Pi=\Pi(H)$ in $H^\circ$. In particular, $H$ has a canonical dual cubing $C(H)$ associated with it -- it is the cubing whose $1$-skeleton is the graph $\Gamma_\Pi$.
	\item Every principal ultrafilter that is not well-founded is the vertex of an infinite-dimensional cube.
\end{enumerate}
\end{prop}
One may get greedy and want to achieve even more, expecting all well-founded ultrafilters to have finite degrees (provided at least one of them has). This is in general untrue: take $H$ to be the poc-set of halfspaces in a bounded tree $T$ with infinitely-many leaves; there is but one almost-equality class, but $T$ inevitably has vertices of infinite degree.

Moreover, it is not even true in general that if $C(H)$ is proper (well-founded ultrafilters exists and are all of finite degree) then $H$ is $\omega$-dimensional: let $T=\{a_n\}_{n\in\NN}\cup\{b_n\}_{n\in\NN}$ be given the partial ordering generated by the relations $b_n\leq a_n$ and $b_n\leq b_{n+1}$ for all $n$, and consider the pocset $H$ generated by $T$. The subset $\{a_n\}_{n\in\NN}$ is an infinite transverse subset, but the cubing $C(H)$ is well-defined and proper, as one may easily verify. Moreover, $H^\circ$ has only one class of principal ultrafilters.\\

The particular geometric realizations of discrete poc-sets that interest us in this paper always have a well-founded ultrafilter, and it is therefore natural to ask if the dual cubing is proper, $\omega$-dimensional, or perhaps even finite dimensional. The above discussion stresses that there is a subtle distinction between the possible $\omega$-dimensionality of $H$ and $C(H)$ not containing infinite-dimensional cubes.

\section{Halfspace systems: definitions.}\label{section:halfspaces}
\begin{figure}[htb]
	\centering{\includegraphics[width=3in]{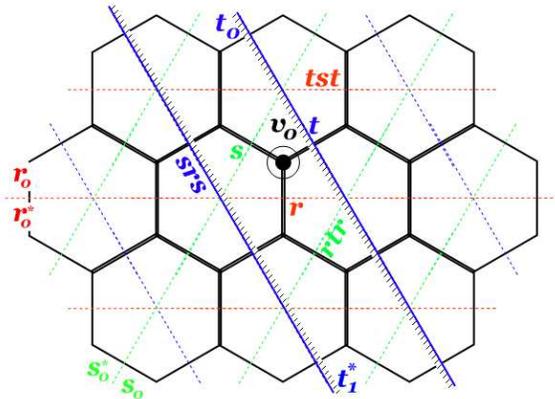}}
	\caption{\protect\scriptsize The hexagonal packing of $\EE^2$: walls are the fixed-point sets of reflections; $\HH$ decomposes as the union of three systems of proper halfspaces $\{r_n\},\{s_n\},\{t_n\}$ (and complements) indexed by $n\in\ZZ$, with $\{r_a,s_b,t_c\}$ transverse for any $a,b,c\in\ZZ$; we set $r_0,s_0,t_0$ to be three pairwise transverse minimal halfspaces among those containing the vertex $v_0$ corresponding to the unit element of $W$.\protect\normalsize}
	\label{figure:halfspace system}
\end{figure}

Our notion of a halfspace system is the obvious generalization of what one observes in a cubing.
\begin{defn}[Halfspace system] A \emph{halfspace system} in a geodesic metric space $X$ is a family $\HH$ of
halfspaces containing the trivial halfspaces, ordered by containment, invariant under the operation $h\mapsto h^\ast$ and satisfying the discreteness condition that every point $x\in X$ has a neighbourhood $U_x$ intersecting only finitely many walls associated with halfspaces of $\HH$.
\end{defn}
A natural example of a CAT(0) space with a halfspace system may be obtained taking $X$ to be the Davis-Moussong complex of a Coxeter system $(W,R)$ of finite rank, and letting $\HH$ be the system of halfspaces arising as the set of complementary components of the walls. A good illustration for most of the work done in this paper is that of the regular hexagonal tiling of the Euclidean plane $\EE^2$, which is nothing else than the Davis-Moussong complex of the Coxeter system
\begin{equation}
	W\cong\left\langle r,s,t\left|r^2,s^2,t^2, (rs)^3,(rt)^3,(st)^3\right.\right\rangle\,,
\end{equation}
as illustrated in figure \ref{figure:halfspace system}.
Walls are defined to be the fixed point sets of reflections of the system $(W,R)$, and it can be shown (for example, see \cite{[Wil]}), that this system of walls coincides with $W(\HH)$.\\

The work of Brink and Howlett \cite{[BriHow]} shows that this particular class of examples has the \emph{parallel walls property}:
\begin{defn}[parallel walls property] A halfspace system $\HH$ in a geodesic metric space $X$ has the \emph{parallel walls property}, if there exists a constant $C>0$ such that for every $h\in\HH$ and $x\in X$ satisfying $d(x,h^\ast)>C$ there exists a halfspace $k\in\HH$ such that $x\in k<h$.
\end{defn}
\begin{remark} When we are mentioning parallel walls, this should not be mistaken for walls lying at a bounded Hausdorff distance from each other. By a pair of parallel walls we only mean walls $W(h),W(k)$ not being separated one from the other by another wall of $\HH$.
\end{remark}
In the applications we have in mind we will need an analogous notion which is formulated with respect to boundary points. This is where CAT(0) geometry sets in:
\begin{defn}[conical points, uniformness] Suppose $\HH$ is a halfspace system in a proper CAT(0) space $X$. A point $\xi\in\bd X$ is said to be a \emph{conical limit point of $\HH$}, if the set
\begin{equation}
	T(\xi)=\left\{h\in\HH\left|\xi\in\textrm{int}\left(\bd h\right)\right.\right\}
\end{equation}
is non-empty, and for any $a\in T(\xi)$ and any cone neighbourhood $U$ of $\xi$ in $X\cup\bd X$ there exists $b\in T(\xi)$ satisfying $b<a$ and $b^\ast\cap U\neq\varnothing$.

A halfspace system $\HH$ on a proper CAT(0) space $X$ is said to be \emph{uniform}, if all points of $\bd X$ are conical limit points of $\HH$.
\end{defn}
Let us verify that a halfspace system with the parallel walls property and satisfying $T(\xi)\neq\varnothing$ for all $\xi\in\bd X$ is uniform.
\begin{lemma}\label{lemma:T(xi) contains descending chains} Suppose $\HH$ is a halfspace system in a proper CAT(0) space $X$, and let $\xi\in\bd X$. If $\HH$ satisfies the parallel walls property, then for every $h\in T(\xi)$ and every cone neighbourhood $U$ of $\xi$ in $\hat{X}$ there exists a $k\in T(\xi)$ such that $k<h$ and $k^\ast\cap U\neq\varnothing$.
\end{lemma}
\proof{} Let $\gamma$ be a geodesic ray in $X$ converging on $\xi$ and emanating from a point $x_0\in h^\ast$. Find $t>0$ such that 
\begin{enumerate}
    \item $d(\gamma(t),h^\ast)>C$, where $C$ is the constant given by the parallel walls property, and such that
    \item the $2C$-neighbourhood of $\gamma([t,\infty))$ is contained in $U$.
\end{enumerate}
By the definition of $C$, there exists $k\in\HH$ satisfying $\gamma(t)\in k<h$. Then, since $\gamma$ crosses $W(k)$ from $k^\ast$ into $k$, we must have $\xi\in\bd k$. Since $\gamma$ is eventually contained in the interior of $k$ (as opposed to the closure of $k^\ast$), we conclude $k\in T(\xi)$ (otherwise, $\gamma$ entering $W(k)$ from $k^\ast$ would have implied $\gamma([t,\infty))\subset W(k)$).

Now, replace $h$ by $k$ and repeat the process if possible (while $d(\gamma(t),h^\ast)>C$). Since only finitely many walls may cross the segment $[x_0,\gamma(t)]$, this process must stop, producing an element $k\in T(\xi)$ containing $\gamma(t)$ and satisfying $d(\gamma(t),k^\ast)\leq C$. For such a $k$, property number (2) of $\gamma(t)$ implies $U\cap k^\ast$ is non-empty.\ep\\

Now we are able to relate uniformness to the situation one encounters for Coxeter groups. Given a halfspace system $\HH$, we notice that most points of the space do not lie on any wall of $\HH$. For any such point $x\in X$ it is possible to associate its {\it (closed) chamber} --
\begin{defn}[generic points, chambers] Let $\HH$ is a halfspace system in a geodesic metric space $X$.
A point $x\in X$ not lying on any wall of $\HH$ is said to be {\it generic} (with respect to $\HH$).
The {\it (closed) chamber} $ch(x)$ of a generic point $x$ is defined as the intersection of closures of all halfspaces in $\HH$ containing $x$.
\end{defn}
In the Davis-Moussong $X=M(W,R)$ complex of a Coxeter system $(W,R)$, all chambers are bounded, as every chamber corresponds to a unique element of $W$, and $W$ acts co-compactly on $X$. We observe --
\begin{prop} Suppose $\HH$ is a halfspace system in a CAT(0) space $X$ such that
\begin{enumerate}
	\item $\HH$ has the parallel walls property, and --
	\item all the chambers of $\HH$ are bounded.
\end{enumerate}
Then $\HH$ is a uniform system.
\end{prop}
\proof{} In view of the preceding lemma it is enough to show that every $\xi\in\bd X$ has $T(\xi)\neq\varnothing$, so suppose $T(\xi)$ is empty for some $\xi$. 

In that case, for every $h\in\HH$ we must have $\xi\in\bd W(h)$. Given a generic point $x$, for any $h\in\HH$ containing $x$ we must then have $[x,\xi)\subset\cl{h}$. This implies $[x,\xi)\subset ch(x)$, contradicting the boundedness of $ch(x)$.\ep

\section{The cubing dual to a halfspace system.}\label{section:cubing dual to a halfspace system} From now on $\HH$ is halfspace system in a geodesic metric space $X$. Since $\HH$ is not assumed to be of dimension $\omega$, the set of principal ultrafilters may split into several almost-equality classes. However, only one of them is directly and naturally associated with the space $X$. The first paragraph in this section resolves this issue and constructs the dual cubing $C(\HH)$, following the original ideas of Sageev. The other two paragraph are new: in the second paragraph we prove Theorem A from the introduction, while the third paragraph develops the technical tools required for proving the finiteness results including Theorem B.

\subsection{Consistent sets and ultrafilters.}\label{subsection:constructing the dual cubing}
The construction in this paragraph follow the scheme laid out by Sageev in \cite{[Sa]} and, more generally, by Nica in \cite{[Nic]} in the setting of discrete spaces.

Henceforth, for any $x\in X$, let $B_x$ denote the set of all $h\in\HH$ containing $x$. Obviously, if $h,k\in B_x$ then $h\not\leq k^\ast$, so that $B_x$ is a filter base in $\HH$ for all $x\in X$. We define:
\begin{defn}\label{defn:consistent set} A point $x\in X$ is said to \emph{support} a subset $A\subset\HH$, if $x\in\cl{h}$ for all $h\in A$.
A non-empty set $A\subset\HH$ is said to be \emph{consistent}, if it has a supporting point.
The empty subset of $\HH$ is, by definition, an inconsistent set.
\end{defn}
\begin{lemma}[consistent ultrafilters]\label{lemma:consistent ultrafilters} Suppose $\HH$ is a halfspace system on a geodesic space $X$.
\begin{enumerate}
	\item A point $x\in X$ supports $\pi\in\HH^\circ$ iff $B_x\subseteq\pi$.
	\item Any consistent ultrafilter is well-founded (and in particular principal).
	\item All consistent ultrafilters lie in the same almost-equality class of $\HH^\circ$.
	\item Any point $x\in X$ supports an ultrafilter.
\end{enumerate}
\end{lemma}
\proof{} To prove (1), Let $\pi\in\HH^\circ$ and $x\in X$. If $x$ supports $\pi$ then every $h\in B_x$ satisfies $x\notin\cl{h^\ast}$, which implies $h^\ast\notin\pi$ and hence $h\in\pi$. Conversely, if $B_x\subseteq\pi$ but $\pi$ is not supported on $x$, then there is $h\in\pi\minus B_x$ such that $x\in h^\ast$; since $h^\ast\in B_x$ we have $h,h^\ast\in\pi$ -- a contradiction.

Note that the above argument amounts to saying that $\pi\minus B_x$ consists of halfspaces $h$ whose walls contain $x$, and is therefore a finite set.

To prove (2), observe that the discreteness requirement on $\HH$ implies the number of $h\in\HH$ whose walls intersect a given compact subset of $X$ must be finite. In particular, no $B_x$ contains an infinite descending chain. By (1), if $\pi\in\HH^\circ$ is consistent then it contains $B_x$ for some $x\in X$. In addition, $\pi\minus B_x$ is finite, thus, since for all $h\in B_x$ the set of $k\in\HH$ satisfying $x\in k<h$ is finite, $\pi$ must be well-founded.

If $\alpha,\beta\in\HH^\circ$ are supported on points $x$ and $y$ respectively, then 
\begin{displaymath}
	\alpha\minus\beta\subseteq(\alpha\minus B_x)\cup(\beta\minus B_y)\cup(B_x\minus B_y)
\end{displaymath}
All three sets on the right hand side are finite, which proves (3).

\begin{figure}[htb]
	\includegraphics[width=\textwidth]{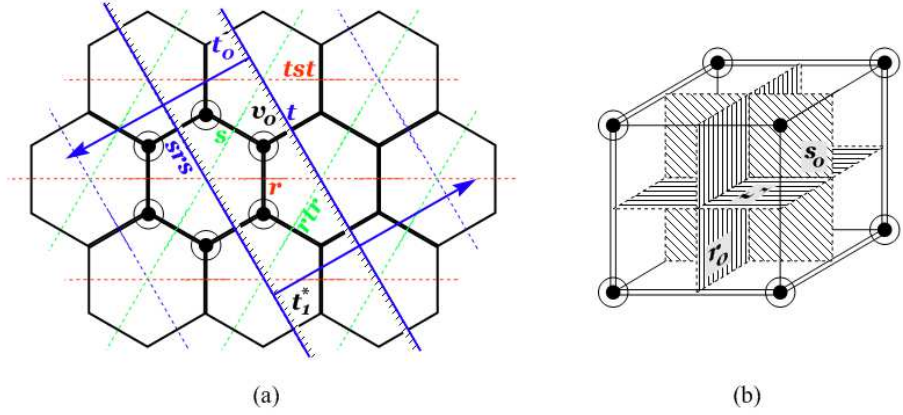}
	\caption{\protect\scriptsize The hexagonal packing of $\EE^2$ (a) demonstrating how inconsistent ultrafilters may arise from transverse sets which are sufficiently `apart' (b).\protect\normalsize}
	\label{figure:inconsistent UF's}
\end{figure}

For any $x\in X$, $B_x$ is a filter base, and is therefore contained in an ultrafilter $\pi$. By (1), $\pi$ is supported on $x$.\ep
\begin{defn}[cubing dual to a halfspace system] Let $\HH$ be a halfspace system in a geodesic metric space $X$. Denote:
\begin{itemize}
	\item[-] by $\Pi_0$ the set of all consistent ultrafilters in $\HH^\circ$;
	\item[-] by $\Pi$ the unique almost-equality class of $\HH^\circ$ that contains $\Pi_0$;
	\item[-] by $C(\HH)$ the cubing whose $1$-skeleton is the graph $\Gamma=\Gamma_\Pi$.
\end{itemize}
In addition, for any $\pi\in\Pi$ we define the {\it height of $\pi$} to be $\hgt{\pi}=\Delta(\pi,\Pi_0)$, and for every $\delta\in\NN\cup\{0\}$ we let $\Pi_\delta$ denote the set of all $\pi\in\Pi$ whose height does not exceed $\delta$.
\end{defn}
Going back to the example of the hexagonal packing in $\EE^2$, figure \ref{figure:inconsistent UF's} shows how a hexagon of the tiling `embeds' in $\Pi$, as well as how inconsistent ultrafilters are formed.

We see that the consistent ultrafilters -- or, at least, those which are of the form $B_x$ for some $x\in X$ -- correspond to the chambers of $\HH$ in $X$. More precisely, the argumentation in lemma \ref{lemma:consistent ultrafilters} shows that a point $x\in X$ is generic (with respect to $\HH$) if and only if $B_x$ is an ultrafilter, and $ch(x)$ is then the set of all points $y$ supporting the ultrafilter $B_x$.

\subsection{Embedding $X$ in the dual cubing.}
Throughout this section, $G$ is a group acting geometrically on a geodesic metric space $X$, and $x_0$ is a fixed basepoint in $X$ which is generic with respect to $\HH$. We also fix a closed ball $B_0=B(x_0,R_0)$ intersecting every $G$-orbit in $X$, and set $\pi_0=B_{x_0}$.

Our goal is to prove the following result:
\begin{prop}[Theorem A]\label{prop:theorem A} If $G$ is a group acting geometrically on a geodesic metric space $(X,d)$, and $\HH$ is a $G$-invariant halfspace system in $X$, then the map
\begin{displaymath}
	\pi:G\cdot x_0\to C(\HH)^1\,,\qquad \pi(x)=B_x
\end{displaymath}
extends to a quasi-isometric embedding of $X$ in $(C(\HH)^1,\Delta)$.
\end{prop}
The idea of the proof is to circumvent the problem of distance estimates using the group $G$ as follows. Let $\Gamma_0$ be the subgraph of $\Gamma=C(\HH)^1$ induced by the vertex set $\Pi_0$, and let $\Delta_0$ be the (combinatorial) path metric on $\Gamma_0$. On the face of it, $\Delta_0\geq\Delta$ and $\Delta_0$ may attain infinite values (if $\Gamma_0$ is disconnected). We shall prove, however, that $\Delta_0$ coincides with $\Delta$. Since, as we shall also see, $\Gamma_0$ is $G$-finite, we can think of $(\Gamma_0,\Delta)$ as a connected geodesic metric space with a co-compact isometric $G$-action, and apply the Milnor schwartz lemma to deduce that the orbit map $\eta:g\mapsto g\cdot\pi_0$ is a quasi-isometry of $G$ with $(\Gamma_0,\Delta)$ (for any fixed finite generating set of $G$). Since the map $\nu:g\mapsto g\cdot x_0$ is a quasi-isometry of $G$ with $X$, the map $\pi$ defined above is a quasi-isometry of $(G\cdot x_0,d)$ with $(\Gamma_0,\Delta)$, which proves the proposition.
\begin{lemma}\label{lemma:consistent UFs are geodesic}
The metric $\Delta_0$ coincides with $\Delta$ on $\Gamma_0$. In particular, $(\Gamma_0,\Delta)$ is a geodesic metric space.
\end{lemma}
\proof{} Let $a,b\in X$. Given ultrafilters $\alpha,\beta$ supported on the points $a$ and $b$ respectively, we construct a geodesic vertex path $(\pi_0=\alpha,\ldots,\pi_k=\beta)$ in $(\Pi,\Delta)$ from $\alpha$ to $\beta$, which is completely contained in $\Pi_0$. The construction is by induction on $\Delta(\alpha,\beta)$.

There is a subdivision $(\left[ p_{i-1},p_i\right])_{i=0}^n$ of $[a,b]$ into subintervals with $p_0=a$ and $p_n=b$ so that none of the open intervals $(p_{i-1},p_i)$ intersects a wall of $H$ separating $a$ from $b$. Consider the set $A=\alpha\minus\beta$, and its decomposition as $A=\bigcup_{i=0}^nA_i$ where $A_i$ consists of those $h\in A$ having $p_i\in W(h)$. Each of the $A_i$ has the ordering induced from $\HH$, and it will be enough to construct an initial segment of the required geodesic path beginning at $\alpha$ and ending at an ultrafilter $\gamma$ satisfying $\Delta(\alpha,\gamma)+\Delta(\gamma,\beta)=\Delta(\alpha,\beta)$ and which is supported on the point $p_1$. 

If $A_0\neq\varnothing$, then in order to obtain $\gamma$ we pick a minimal element $h_1$ of $A_0=A_0^{(0)}$, and set $\pi_1=[\pi_0]_{h_1}$, $A_0^{(1)}=A_0\minus\{h_1\}$, and for any $1\leq t\leq |A_0|$, we pick a minimal element $h_t\in A_0^{(t-1)}$ and set $\pi_t=[\pi_{t-1}]_{h_t}$ and $A_0^{(t)}=A_0^{(t-1)}\minus\{h_t\}$. The resulting path from $\alpha=\pi_0$ to $\gamma=\pi_{|A_0|}$ is, by construction, a geodesic vertex path with $\gamma$ satisfying both the required properties.  

If $A_0=\varnothing$, then $\alpha$ is supported on $p_1$, and we may apply the same procedure to the interval $[p_1,b]$ and the set $A_1$.\ep
\begin{lemma}\label{lemma:consistent principal UF's are G-finite}
The set $\Pi_0$ is $G$-finite.
\end{lemma}
\proof{} Let $K$ be a compact ball intersecting every $G$-orbit in $X$, and let $A_K$ be the set of all $h\in\HH$ whose walls intersect $K$. It will be enough to show that the set of all principal ultrafilters supported on $K$ is finite. Let $B_K$ be the set of all $h\in\HH$ containing $K$. Then $B_K$ is contained in any element of $\Pi_0$ that is supported on $K$,  implying that any two $\sigma,\sigma'\in\Pi_0$ containing $B_K$ satisfy $\sigma\vartriangle\sigma'\subseteq A_K$. Since $A_K$ is a finite set, we are done.\ep\\

This concludes the proof of proposition \ref{prop:theorem A}.

\subsection{Studying the height function.}
Our motivation for studying the height function comes from the problem of characterizing the situations when the quotient of $C(\HH)$ by the action of $G$ is compact. Of course, this happens iff $G$ acts co-finitely on $\Pi=C(\HH)^0$. Clearly, the metric $\Delta$ on $\Pi$ is $G$-invariant, and since $G$ acts on $\Pi$ stabilizing $\Pi_0$, $G$ also stabilizes the sub-level sets $\Pi_\delta$ ($\delta\in\NN\cup\{0\}$) of the height function. The action of $G$ on $\Pi_0$ is co-finite, so the action of $G$ on $\Pi$ will be co-bounded if and only if the height function is bounded. Thus, studying the growth of the height function is necessary for studying any finiteness properties (local or global) that the complex $C(\HH)$ may have.

\subsubsection{Distance to $\Pi_0$.}
We need a tool for computing the height function.
\begin{lemma}\label{lemma:characterizing level sets} Suppose $\pi\in\Pi$. Then $\hgt{\pi}\leq n$ if and only if there exists a subset $A$ of $\pi$ of size $n$ such that $\pi\minus A$ is consistent.
\end{lemma}
\proof{} Suppose $A\subseteq\pi$ admits a point $x\in X$ supporting $\pi\minus A$.

For any $h\in B_x$ we have $x\in h$. If $h\notin\pi$, then $h^\ast\in\pi$ and there are two cases to consider:
\begin{itemize}
    \item $h^\ast\in\pi\minus A$. This implies $x\in\cl{h^\ast}$ -- a contradiction.
    \item $h^\ast\in A$. This is the same as $h\in A^\ast$.
\end{itemize}
Thus, $B_x\minus\pi$ lies in $A^\ast$.

Next, let $\pi_x$ be an ultrafilter which is supported on $x$ and such that $\Delta(\pi_x,\pi)$ is minimal. We may write:
\begin{eqnarray}
    \pi_x\minus\pi&=&(B_x\minus\pi)\cup((\pi_x\minus B_x)\minus\pi)\\
    &\subseteq& A^\ast\cup\underbrace{(\pi_x\minus B_x)\minus\pi}_{(\ast)}.
\end{eqnarray}
We will show that $(\ast)$ is the empty set. If not, then select a minimal element $h$ of $(\ast)$.

We first claim $h\in\min(\pi_x)$. For suppose $k\in\pi_x$ satisfies $k<h$: $k<h$ is the same as $h^\ast<k^\ast$, and since $h^\ast\in\pi$, we also have $k^\ast\in\pi$, so that $k\in\pi_x\minus\pi$; by the minimality assumption regarding $h$, this may happen only in case $k\in B_x$; thus, on one hand we have that $k^\ast\in\pi$ forces $x\in\cl{k^\ast}$, while on the other hand we obtain $x\in k$ -- a contradiction.

Thus, the ultrafilter $[\pi_x]_h$ exists. By construction, $[\pi_x]_h$ is an ultrafilter supported on $x$, whose distance to $\pi$ is by one smaller than the allowed minimum -- again, a contradiction, -- and we conclude $(\ast)$ must be the empty set. The inclusion we have consequently obtained shows then that $\Delta(\pi_x,\pi)\leq|A|$, as desired.\\

Conversely, given $\pi\in\Pi$, suppose there is a point $x\in X\minus\bigcup_{h\in\HH}W(h)$ such that $\Delta(\pi,\pi_x)\leq n$ for some ultrafilter $\pi_x$ containing $B_x$. Setting $A$ to be any subset of $\pi$ of size $n$ and containing $\pi\minus\pi_x$ will result in $\pi\minus A$ being a consistent set (because it is a subset of $\pi_x$, which is consistent).\ep

\subsubsection{Height growth and Shadows.}
Let us study how $\hgt{\pi}$ changes as $\pi$ ``moves around'' $\Pi$. For this we need some technicalities.
\begin{defn} Suppose $\pi\in\Pi$ and $a\in\min(\pi)$. Denote
\begin{eqnarray*}
	a\in\min(\pi)_+ &\IFF& \hgt{[\pi]_a}>\hgt{\pi}\,,\\
	a\in\min(\pi)_- &\IFF& \hgt{[\pi]_a}<\hgt{\pi}\,,\\
	a\in\min(\pi)_0 &\IFF& \hgt{[\pi]_a}=\hgt{\pi}.
\end{eqnarray*}
Note that $\pi\in\Pi_0$ iff $\min(\pi)_-$ is empty.
\end{defn}
\begin{lemma}\label{lemma:min pi minus is inconsistent} For all $\pi\in\Pi$, the set $\min(\pi)_-$ is inconsistent.
\end{lemma}
\begin{cor}\label{cor:going down transversely}
Suppose $\pi\in\Pi\minus\Pi_0$, and $a\in\pi$. Then there exists $c\in\min(\pi)_-$ such that $c\pitchfork a$.
\end{cor}
\proof{} if all $b\in\min(\pi)_-$ were facing $a$, we would have $a^\ast<b$ for all such $b$, contradicting the lemma.\ep

\proof{of lemma \ref{lemma:min pi minus is inconsistent}} By induction on $\delta=\hgt{\pi}$: for $\delta=0$ the statement is trivial, so assume $\pi\notin\Pi_0$ and that any ultrafilter at a distance $\delta-1$ satisfies the statement of this lemma.

Now, by lemma \ref{lemma:characterizing level sets}, there exists an $a\in\min(\pi)_-$. We consider elements $b\in\min([\pi]_a)_-$: for every such $b$ we must either have $a<b$ or $b\in\min(\pi)_-$; as a result we obtain the containment
\begin{equation}
	\bigcap_{b\in\min(\pi)_-}\cl{b}\subseteq\bigcap_{b\in\min([\pi]_a)_-}\cl{b}.
\end{equation} 
Since the right-hand side is empty (induction hypothesis), so is the left-hand side.\ep

\begin{cor}\label{cor:three ways to go down} If $\pi\notin\Pi_0$, then $\min(\pi)_-$ contains at least three distinct elements.
\end{cor}
\proof{} This is immediate: $\min(\pi)_-$ is, first of all, a filter base, and hence every pair of elements in $\min(\pi)_-$ is consistent; therefore, in order to be inconsistent, it must contain at least three elements.\ep

\begin{defn}[Shadows]\label{defn:shadow} For all $\pi\in\Pi$, let the \emph{shadow of $\pi$} be defined as the set
\begin{equation}
	\sh{\pi}=\left\{\sigma\in\Pi_0\,\big|\,\Delta(\sigma,\pi)=\hgt{\sigma}\right\},
\end{equation} 
and let the \emph{dual shadow} $\shcirc{\pi}$ of $\pi$ be defined to be
\begin{equation}
	\shcirc{\pi}=\left\{h\in\HH\,\big|\,\sh{\pi}\subseteq S_h\right\}.
\end{equation}
\end{defn}
\begin{figure}[htb]
	\includegraphics[width=\textwidth]{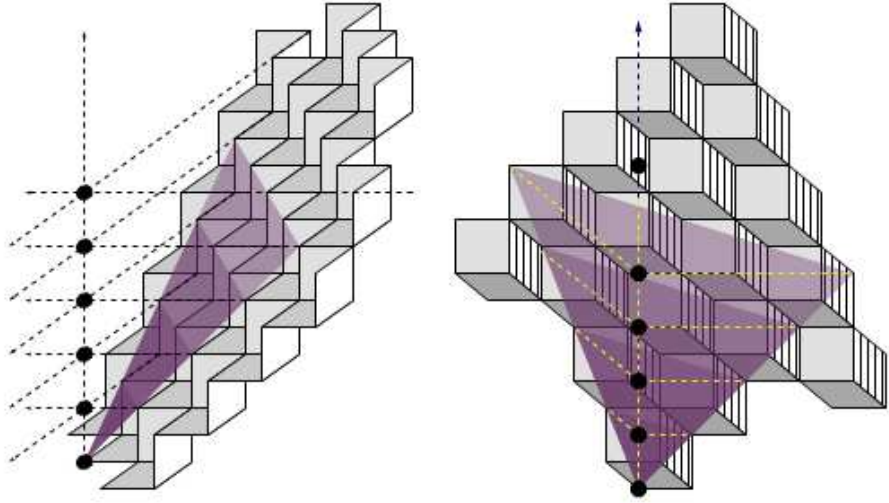}
	\caption{\protect\scriptsize the `shadows' cast by ultrafilters at increasing distance to $\Pi_0$, for the case of the cubing dual to the hexagonal tiling of $\EE^2$. The painted $2$-cubes are the ones whose vertices lie in $\Pi_0$.\protect\normalsize}
	\label{figure:shadows grow}
\end{figure}
Observe that $\shcirc{\pi}$ is a filter-base, and it is natural to expect that $\shcirc{\pi}$ be contained in $\pi$. If that is the case, it would mean that, as the distance of $\pi$ from $\Pi_0$ increases, $\shcirc{\pi}$ diminishes accordingly, testifying to the growth of $\sh{\pi}$. This also provides one with a tool to assess the size and positioning of the shadow `cast' by a given ultrafilter $\pi$. 

The next few observations are motivated by the example of how $\Pi_0$ embeds in $\Pi$ in the case of the hexagonal packing in $\EE^2$ -- see figure \ref{figure:shadows grow}. They become an important technical tool in what follows.
\begin{lemma}\label{lemma:shadows grow} Suppose $\pi\in\Pi$ and $h\in\min(\pi)_+$. Then $\sh{\pi}\subseteq\sh{[\pi]_h}$ and, consequently, $\shcirc{[\pi]_h}\subseteq\shcirc{\pi}$.
\end{lemma}
\proof{} Observe that the inclusion $\sh{\pi}\subseteq\sh{[\pi]_h}$ automatically implies the reverse inclusion of the dual shadows.

Now, since changing the orientation of $h$ in $\pi$ increases the distance of $\pi$ to $\Pi_0$ by $1$ and changes the distance of $\pi$ to any other $\sigma\in\Pi$ by exactly $1$, it follows that for any $\sigma\in\sh{\pi}$ we must have $\hgt{[\pi]_h}=\hgt{\pi}+1$. In particular, $\sigma\in\sh{[\pi]_h}$.\ep
\begin{lemma}\label{lemma:going down avoiding inversions} For every $\pi\in\Pi$ and every $a,b\in\pi$ there exists $\sigma\in\sh{\pi}$ satisfying $a,b\in\sigma$.
\end{lemma}
\proof{} For fixed $a,b\in\HH$ and $\pi\in\Pi$ containing the halfspaces $a$ and $b$ there exists $c\in\min(\pi)_-$ such that $c\neq a,b$. Thus, $a,b\in[\pi]_c$ and the distance of $\pi$ to $\Pi_0$ is reduced by $1$. Repeating this procedure $\hgt{\pi}$ times we obtain the required $\sigma$.\ep
\begin{lemma}\label{lemma:bounds on shadows} For every $\pi\in\Pi$ one has
$\min(\pi)_+\cup\min(\pi)_0\subseteq\shcirc{\pi}\subseteq\pi$.
\end{lemma}
\proof{} For the trivial case when $\pi\in\Pi_0$ one has $\sh{\pi}=\{\pi\}$, $\shcirc{\pi}=\pi$ and $\min(\pi)_-$ is empty, proving the required inclusions.

When $\pi\notin\Pi_0$, let us verify the left-hand side first. If $\sigma\in\sh{\pi}$ and $a^\ast\in\sigma$, then we must have
\begin{equation}
	\Delta([\pi]_a,\sigma)=\Delta(\pi,\sigma)-1=\hgt{\pi}-1,
\end{equation}
implying $a\in\min(\pi)_-$. This proves the implication $a\notin\shcirc{\pi}\THEN a\notin\min(\pi)_+\cup\min(\pi)_0$.\\

For the right hand side inclusion we take $k\notin\pi$ and show that $k\notin\shcirc{\pi}$. Since $k^\ast\in\pi$, the last lemma provides $\sigma\in\sh{\pi}$ containing $k^\ast$, which proves $k\notin\shcirc{\pi}$.\ep\\

In the light of the last result, consider $\pi\in\Pi$ and $a\in\min(\pi)_+$: we then know that $a\in\shcirc{\pi}$, while $a^\ast\in[\pi]_a$ and the last lemma imply that $a\notin\shcirc{[\pi]_a}$. Also, $a^\ast\notin\shcirc{[\pi]_a}$, because $\shcirc{[\pi]_a}$ is contained in $\shcirc{\pi}$, which is a filter base containing $a$. In particular, this means that $\sh{\pi}\subsetneq\sh{[\pi]_a}$. As a result, we obtain that the shadow of $\pi$ strictly increases as $\pi$ is moved farther and farther away from $\Pi_0$. Observe, in the above discussion, that $a\in\min([\pi]_a)_-$, so this result could also be stated by saying that $\shcirc{\pi}$ is disjoint from $\min(\pi)_-$.\\

The above technical observations are summarized in the following proposition --
\begin{cor}[strict growth of shadows]\label{cor:shadows really grow} Suppose $\pi\in\Pi$. Then,
\begin{enumerate}
	\item if $a\in\min(\pi)_+$, then $\sh{\pi}\subsetneq\sh{[\pi]_a}$ and $\shcirc{[\pi]_a}\subsetneq\shcirc{\pi}$;
	\item $\shcirc{\pi}\cap\min(\pi)=\min(\pi)_+\cup\min(\pi)_0$.\,\ep
\end{enumerate}
\end{cor}

\section{The geometry of uniform systems.}\label{section:geometry of uniform systems}
Hereafter $G$ is a group acting geometrically on a CAT(0) space $X$ (which is therefore proper), and $\HH$ is a uniform $G$-invariant halfspace system.

We seek some understanding of the finiteness properties of the cubing $C(\HH)$ dual to $\HH$, as well as of $\HH$ itself. Once again, we let $\Gamma$ denote the $1$-skeleton of this cubing, and recall that the vertex set of $\Gamma$ is $V\Gamma=\Pi$ and the metric $\Delta$ defined on $\Pi$ is precisely the combinatorial metric on $\Gamma$ assigning unit length to all edges.

We now consider a geometric counterpart of shadows.
\begin{defn}[geometric shadow] For $\pi\in\Pi$, define the \emph{geometric shadow} $\gsh{\pi}$ of $\pi$ to be the support of $\shcirc{\pi}$:
\begin{equation}
	\gsh{\pi}=\bigcap_{a\in\shcirc{\pi}}\cl{a}\,.
\end{equation}
\end{defn}
Clearly, if $x$ supports an element of $\sh{\pi}$, then $x\in\gsh{\pi}$. Combined with the technical results we have on shadows, this notion serves as the main tool for relating the geometry of $C(\HH)$ to that of $X$. 

\subsection{The main theorem.}
We are now able to prove our main theorem.
\begin{thm}[theorem B]\label{thm:dual cubing is locally finite} Suppose $G$ is a group acting geometrically on a proper CAT(0) space $X$, and suppose $\HH$ is a uniform $G$-invariant halfspace system in $X$. Then the action of $G$ on $\Pi_\delta$ is co-finite for every $\delta\in\NN$ and the cubing $C(\HH)$ is locally-finite. In particular, $C(\HH)$ contains no infinite-dimensional cube, and the action of $G$ on $C(\HH)$ is co-compact if and only if the height function is bounded.
\end{thm}
\proof{} Let $x_0\in X$ be a fixed base point and let $R_0>0$ be such that the closed ball $B_0=B_d(x_0,R_0)\subset X$ intersects every orbit of $G$ in $X$. For each $\delta\in\NN$ let us denote the set of ultrafilters $\pi\in\Pi$ with $\Delta(\pi,\Pi_0)\leq\delta$ and $\gsh{\pi}\cap B_0\neq\varnothing$ by $S(\delta)$.

The first step of the proof involves producing a bound on the diameter of geometric shadows: we claim that for every $\delta\in\NN$ there exists $R(\delta)>0$ such that for every $\pi\in\Pi$ satisfying $\Delta(\pi,\Pi_0)\leq\delta$ there exists $g\in G$ satisfying
\begin{equation}
	\gsh{g\cdot\pi}=g\cdot\gsh{\pi}\subseteq B\left(x_0,R(\delta)\right)\,.
\end{equation}
Fix $\delta\in\NN$. Since $B_0$ intersects every orbit of $G$, it is enough to prove there exists a number $R>0$ satisfying $\gsh{\pi}\subseteq B\left(x_0,R\right)$ for any $\pi\in S(\delta)$.

Let $\pi\in S(\delta)$ and consider a point $x\in\gsh{\pi}$ such that $d(x,x_0)\leq R_0$, and a boundary point $\xi\in\bd X$. By uniformness, there exists a descending sequence $a_1(\xi),\ldots,a_N(\xi)\in T(\xi)$ such that $B_0\subset a_1^\ast$ and $N=2\delta+1$.

We claim $a_N(\xi)\cap\gsh{\pi}$ is empty. For suppose $y\in a_N(\xi)\cap\gsh{\pi}$. Then $x\in a_1(\xi)^\ast$ and $y\in a_1(\xi)$ implies neither $a_1(\xi)$ nor $a_1(\xi)^\ast$ lie in $\shcirc{\pi}$. In the same manner we conclude that $a_N(\xi),a_N(\xi)^\ast$ do not lie in $\shcirc{\pi}$. Then, there exist $\alpha,\beta\in\sh{\pi}$ with $a_1(\xi)^\ast\in\alpha$ and $a_N(\xi)\in\beta$. But then $a_i(\xi)^\ast\in\alpha$ and $a_i(\xi)\in\beta$ for all $i=1,\ldots,N$, and we conclude that 
\begin{equation}
	\Delta(\alpha,\beta)\geq N>2\delta\,.
\end{equation}
However, this is impossible, as
\begin{equation}
	\Delta(\alpha,\beta)\leq\Delta(\alpha,\pi)+\Delta(\pi,\beta)=2\delta\,.
\end{equation} 

Now, fixing a sequence $a_1(\xi),\ldots,a_N(\xi)$ as above for every $\xi\in\bd X$, recall that each $a_N(\xi)$ contains a cone neighbourhood of $\xi$ (in $X$). Since $X\cup\bd X$ is compact in the cone topology ($X$ is proper), there exists $R>0$ such that the ball $B_d(x_0,R)$ contains $X\minus\bigcup_{\xi\in\bd X}a_N(\xi)$. Thus, the preceding calculation shows that $y\notin\gsh{\pi}$ whenever $y\notin B_d(x_0,R)$, and the first step is done.\\

The second step of the proof relates balls in $\Pi$ (defined by the metric $\Delta$) to balls in $X$. Again, we fix some natural number $\delta$. Consider now an ultrafilter $\pi\in B_\Delta\left(\sigma,\delta\right)$, where $\sigma$ is an ultrafilter supported on a point $x\in X$. Let $y$ be a point supporting an element $\pi_y\in\sh{\pi}$ (and so, $y\in\gsh{\pi}$), and find $g\in G$ such that $g\cdot y\in B_0$. By the construction of $R(\delta)$, if $d(x,y)>R(\delta)$ then $\Delta(g\cdot\sigma,g\cdot\pi_y)$ is greater than or equal to $2\delta+1$. On the other hand, we have
\begin{equation}
	\Delta(g\cdot\sigma,g\cdot\pi_y)\leq
	\Delta(g\cdot\sigma,g\cdot\pi)+\Delta(g\cdot\pi,g\cdot\pi_y)=
	\Delta(\sigma,\pi)+\Delta(\pi,\pi_y)\leq
	2\delta,
\end{equation} 
producing a contradiction again. In particular, if $x=x_0$ (is the point supporting $\sigma$) and $\Delta(\pi,\sigma)\leq\delta$ then $\gsh{\pi}$ intersects $B_d\left(x_0,R(\delta)\right)$; however, since $\Delta(\pi,\Pi_0)\leq\delta$ we also know that the diameter of $\gsh{\pi}$ is at most $2R(\delta)$, which implies $\gsh{\pi}$ is contained in the closed ball $B_d\left(x_0,3R(\delta)\right)$. This concludes the second step.\\

Now, For each $\pi\in S(\delta)$ recall that $\shcirc{\pi}$ is contained in $\pi$, and consider $h\in\pi\minus\shcirc{\pi}$. For such an $h$, we must have both $S_h\cap\sh{\pi}$ and $S_{h^\ast}\cap\sh{\pi}$ non-empty, providing us with ultrafilters $\sigma,\sigma^\ast\in\sh{\pi}$ satisfying $h\in\sigma$ and $h^\ast\in\sigma^\ast$. We conclude that there exist points $x\in h$ and $x\in h^\ast$, both lying in the \emph{open} ball $B\left(x_0,R(\delta)+1\right)$. Thus, $\pi\minus\shcirc{\pi}$ is contained in the subset of all halfspaces of $\HH$ whose walls intersect a ball about $x_0$ whose radius depends only on $\delta$, implying that $\pi\minus\shcirc{\pi}$ is a finite set of size bounded by a function of $\delta$. Now, since $\gsh{\pi}$ is contained in that same ball, we see that there also are only finitely many possibilities for selecting $\shcirc{\pi}$ given $\delta$. We have shown that $S(\delta)$ is a finite set, and the first step then allows the conclusion that $G$ acts co-finitely on the level set $\Pi_\delta$ of the function $\Delta(-,\Pi_0)$.

Finally, using the second step, let us employ the properness of the action of $G$ on $X$ to deduce local finiteness. Consider the ball $B_\Delta\left(\sigma_0,\delta\right)$ in $\Pi$, where $\sigma_0$ is an ultrafilter supported on the basepoint $x_0$. If $(\pi_n)_{n=1}^\infty$ are all distinct elements of $B_\Delta(\sigma_0,\delta)$, and $g_n\in G$ are such that $g_n\cdot\pi_n\in S(\delta)$ for all $n$, then the finiteness of $S(\delta)$ allows passing to a subsequence in which $g_n\cdot\pi_n=\pi$ for all $n$, for some suitable $\pi\in S(\delta)$; considering the ball $B=B_d\left(x_0,3R(\delta)\right)$, we then have (by the second step) for all $n$
\begin{equation}
	g_n^{-1}g_1\cdot B\supseteq g_n^{-1}g_1\cdot\gsh{\pi_1}=\gsh{\pi_n}\subseteq B,
\end{equation}
contradicting the proper-discontinuity of the action of $G$ on $X$.\ep

\subsection{Parallel walls vs. uniformness.}
For every $r\geq 0$, denote the closed ball of radius $r$ about $x_0$ by $B(r)$, and set $B_0=B(R_0)$, where, as before, $R_0$ is chosen so that $B(R_0)$ intersects every orbit of $G$ in $X$. Further let $f(r)$ denote the minimum, over all $k\in\HH$ containing $B(r)$, of the number of $h\in\HH$ containing $B_0$ and satisfying $h\leq k$, and set $T_0$ to be the number of walls of $\HH$ intersecting $B_0$.

Let now $x\in X$ and $r>0$, and let $k\in\HH$ contain $B(x,r+R_0)$. Find $g\in G$ such that $g\cdot x\in B_0$, so that $g\cdot k$ contains $B(r)$. Then, every $h\in\HH$ containing $B_0$ and satisfying $h\leq g\cdot k$ will satisfy $x\in g^{-1}\cdot h\leq k$. Thus, the set $S(x,k)$ of all $h\in\HH$ satisfying $x\in h\leq k$ contains at least $f(r)$ distinct elements. We summarize this in 
\begin{lemma}[counting walls]\label{lemma:counting walls} There exists a non-decreasing function $f:\RR_+\to\NN$ such that every $r>0$, $x\in X$ and $k\in\HH$ satisfy
\begin{equation}
	d(x,k^\ast)>r+R_0\quad\THEN\quad \left|S(x,k)\right|\geq f(r),
\end{equation}
where $S(x,k)$ is the set of all $h\in\HH$ satisfying $x\in h\leq k$.
\end{lemma}
The function $f$ reflects various finiteness properties of $\HH$. For example, the parallel walls property for the pair $(X,\HH)$ is equivalent to saying that $f(r)$ is strictly greater than $1$ for a big enough value of $r$. Thus, we are interested in information regarding the growth of the function $f$. To that end, we introduce the following notion:
\begin{defn}[slope of a halfspace system] Suppose $G$ is a group acting geometrically on a CAT(0) space $X$, and $\HH$ is a $G$-invariant uniform halfspace system on $X$. We say that $\HH$ has slope $\alpha\geq 0$ if $\liminf_{r\to\infty}{f(r)/r}=\alpha$.
\end{defn}
Note that negating the parallel walls property implies $\HH$ has zero slope.\\

Let us consider the set $S(x_0,k)$ for some $k$ containing $x_0$. For any $\sigma\in S_{k^\ast}$, it is clear that $S(x_0,k)\subseteq \pi_0\minus\sigma$, implying $\Delta\left(\pi_0,S_{k^\ast}\right)$ is at least $|S(x_0,k)|$. In the reverse direction, if $\sigma_0$ is the projection of $\pi_0$ to $S_{k^\ast}$, let us compute $\pi_0\minus\sigma_0$:
\begin{itemize}
	\item[-] If $a\in\pi_0$ satisfies $a\pitchfork k$, then we will also have $a\in\sigma_0$. Indeed, recall that for any $\alpha\in S_{k^\ast}$ we have $\sigma_0\in\left[\pi_0,\alpha\right]$; since $a\pitchfork k$, we may choose $\alpha\in V(a,k^\ast)$, which then implies $\sigma_0\in V(a,k^\ast)$.
	\item[-] If $a\in\pi_0$ satisfies $a^\ast<k$, then $S_{k^\ast}\subset S_a$, so that $\sigma_0\in S_a$.
	\item[-] No $a\in\pi_0$ may satisfy $a<k^\ast$.
	\item[-] If $a\in\pi_0$ satisfies $k<a$, then $\sigma_0$ will contain $a$. Indeed, consider $\sigma_1=[\sigma_0]_{k^\ast}$: this is an ultrafilter containing $k$, and so contains $a$, which implies that $[\sigma_1]_{k}=\sigma_0$ contains $a$, too.
	\item[-] Finally, we conclude from the above that the only elements of $\HH$ possibly contained in $\pi_0\minus\sigma_0$ are the elements of $S(x_0,k)$.
\end{itemize}
We have proved
\begin{lemma}[distance to a wall]\label{lemma:distance to a wall} for any generic point $x\in X$ and any $k\in\pi_x$ one has 
\begin{equation}
	\Delta\left(\pi_x,S_{k^\ast}\right)=\left|S(x,k)\right|\,.\quad\ep
\end{equation}
\end{lemma}
\begin{cor}[uniform implies parallel walls property] Suppose $G$ is a group acting geometrically on a CAT(0) space $X$, and $\HH$ is a $G$-invariant uniform halfspace system on $X$. Then $\HH$ has the parallel walls property.
\end{cor}
\proof{} In view of lemma \ref{lemma:counting walls} and the fact that $\left||S(x_0,k)|-f(d(x_0,k^\ast))\right|\leq T_0$, the negation of the parallel walls property is equivalent to having $\{k_n\}_{n\in\NN}\subset\HH$ with the property that $S(x_0,k_n)=\{k_n\}$ for all $n$ and $d(x_0,k_n^\ast)\to \infty$ as $n\to\infty$. 

For each $n$, consider the ultrafilter $\pi_n=[\pi_0]_{k_n}$. This is an inconsistent ultrafilter, and therefore lies in $\Pi_1$. By lemma \ref{lemma:going down avoiding inversions}, there exists $a_n\in\min(\pi_n)_-$ such that $a_n\neq k_n^\ast$. Thus, if $x_n$ is a point of $X$ supporting $\sigma_n=[\pi_n]_{a_n}$, then $x_n\in\cl{k_n^\ast}$. Since $x_n\in\gsh{\sigma_n}$ and $x_0\in\gsh{\pi_0}$, we have that both $x_0$ and $x_n$ lie in $\gsh{\pi_n}$. Therefore, $(\pi_n)_{n=1}^\infty$ is a sequence of ultrafilters in $\Pi_1$ having geometric shadows of unbounded diameter -- a contradiction to the first stage of the proof of theorem \ref{thm:dual cubing is locally finite}.\ep\\

The argument above can be generalized. Suppose $k\in\pi_0$ and $\delta\geq 0$ are such that the projection $\pi$ of $\pi_0$ to $S_{k^\ast}$ satisfies $\Delta(\pi,\Pi_0)=\delta$. Using lemma \ref{lemma:going down avoiding inversions} to avoid inverting $k^\ast$, we can construct a path of length $\delta$ from $\pi$ to $\Pi_0\cap S_{k^\ast}$. As a result we have:
\begin{equation}\label{eqn:lower bound on slope1}
	\Delta\left(\pi_0,S_{k^\ast}\right)\leq\Delta\left(\pi_0,\Pi_0\cap S_{k^\ast}\right)\leq\Delta(\pi_0,\pi)+\delta
		\leq 2\Delta\left(\pi_0,S_{k^\ast}\right)
\end{equation}
Now, let $x_1$ be the projection of $x_0$ to the closure of $k^\ast$. There is a constant $\epsilon'$ independent of $k$ such that $\pi(x_1)$ is at a distance $\epsilon'$ from an element of $\Pi_0\cap S_{k^\ast}$ which is also supported on $x_1$. Since $\pi$ is a $(\lambda,\epsilon)$-quasi isometry of $(X,d)$ with $(\Pi_0,\Delta)$, we also have that the following inequalities hold:
\begin{eqnarray*}
	\lambda\cdot d(x_0,x_1)
		&\geq&\Delta\left(\pi(x_0),\pi(x_1)\right)-\epsilon\\
		&\geq&\Delta\left(\pi_0,\pi_1\right)-\epsilon'-\epsilon\\
		&\geq&\Delta\left(\pi_0,\Pi_0\cap S_{k^\ast}\right)-D-\epsilon\,.
\end{eqnarray*}
This implies
\begin{equation}
	\limsup_{d(x_0,k^\ast)\to\infty}\frac{\Delta\left(\pi_0,\Pi_0\cap S_{k^\ast}\right)}{d(x_0,k^\ast)}\leq\lambda\,.
\end{equation}
On the other hand, for every $x\in k^\ast$ we have
\begin{eqnarray*}
	\Delta\left(\pi(x_0),\pi(x)\right)&\geq&\frac{1}{\lambda}d(x_0,x)-\epsilon
		\geq\frac{1}{\lambda}d(x_0,x_1)-\epsilon,
\end{eqnarray*}
producing the inequality
\begin{equation*}
	\Delta\left(\pi_0,\Pi_0\cap S_{k^\ast}\right)+\epsilon'\geq\frac{1}{\lambda}d(x_0,k^\ast)-\epsilon
\end{equation*}
In the limit we shall then have
\begin{equation}
	\liminf_{d(x_0,k^\ast)\to\infty}\frac{\Delta\left(\pi_0,\Pi_0\cap S_{k^\ast}\right)}{d(x_0,k^\ast)}\geq\frac{1}{\lambda}\,.
\end{equation}
Thus, for the slope of $\HH$ we may then write down the inequalities:
\begin{eqnarray*}
	\liminf_{r\to\infty}\frac{f(r)}{r}&\geq&
	\lambda^{-1}\liminf_{d(x_0,k^\ast)\to\infty}
		\frac{\Delta\left(\pi_0,S_{k^\ast}\right)}{\Delta\left(\pi_0,\Pi_0\cap S_{k^\ast}\right)}\\
	\limsup_{r\to\infty}\frac{f(r)}{r}&\leq&
	\lambda\limsup_{d(x_0,k^\ast)\to\infty}
		\frac{\Delta\left(\pi_0,S_{k^\ast}\right)}{\Delta\left(\pi_0,\Pi_0\cap S_{k^\ast}\right)}\,.
\end{eqnarray*}
However, by equation \ref{eqn:lower bound on slope1}, this implies 
\begin{equation}\label{eqn:bounds on slope}
	\frac{1}{2\lambda}\leq
	\liminf_{r\to\infty}\frac{f(r)}{r}\leq\limsup_{r\to\infty}\frac{f(r)}{r}\leq
	\lambda\,.
\end{equation}
We have proved:
\begin{prop} If $G$ is a group acting geometrically on a CAT(0) space $X$, and $\HH$ is a $G$-invariant uniform halfspace system on $X$, then $f(r)$ has linear growth. In particular, $\HH$ has positive slope.
\end{prop}
This means that $\left|S(x,k)\right|$ grows linearly as a function of $d(x,k^\ast)$ with the rate independent of $x$, which we think may serve as an indication to $\HH$ satisfying an even stronger parallel walls property:
\begin{defn}[\cite{[NibRee1]}] A halfspace system $\HH$ is said to have the \emph{strong} parallel walls property, if there exists a constant $K$ such that for all $a,b\in\HH$ satisfying $d(b^\ast,a)>K$ there exists $h\in\HH$ such that $a<h<b$.
\end{defn}
For the purpose of our discussion of this property we introduce a technical term
\begin{defn} A pair of halfspaces $a,b\in\HH$ is said to be {\it special}, if $a<b$ and there is no $h\in\HH$ satisfying $a<h<b$. The width of such a pair is defined to be $w(a,b)=d(a,b^\ast)$.
\end{defn}
Let us consider a special pair $(a,b)$. By lemma 2.14 in \cite{[Rol]}, for any $a,b\in\HH$ there exist ultrafilters $\alpha,\beta^\ast\in\Pi$ satisfying
\begin{itemize}
	\item[-] $\alpha\in S_a$, $\beta^\ast\in S_{b^\ast}$;
	\item[-] $\Delta(\alpha,S_{b^\ast})=\Delta(\beta,S_a)=\Delta(S_a,S_{b^\ast})$;
	\item[-] $\Delta(\alpha,\beta^\ast)=\left|\left\{h\in H\left| S_a\subseteq S_h\subseteq S_b\right.\right\}\right|$
\end{itemize}

In the terminology of \cite{[Rol]}, lemma 2.13, the pair $(\alpha,\beta^\ast)$ is a {\it gate} for the pair $(S_a,S_{b^\ast})$ of (convex) subsets of the median algebra $\Pi$. In the case of a special pair we will have $\Delta(\alpha,\beta^\ast)=2$. Suppose now that $\delta\in\NN$ is such that $\Delta(\alpha,\Pi_0),\Delta(\beta^\ast,\Pi_0)\leq\delta$. By lemma \ref{lemma:going down avoiding inversions}, we know there exist $\alpha_0\in\sh{\alpha}\cap S_a$ and $\beta_0\in\sh{\beta^\ast}\cap S_{b^\ast}$ and so it follows that:
\begin{eqnarray*}
	\frac{1}{\lambda}w(a,b)-\epsilon&\leq&\Delta(\alpha_0,\beta_0)\\
		&\leq& \Delta(\alpha_0,\alpha)+\Delta(\beta_0,\beta^\ast)+\Delta(\alpha,\beta^\ast)\\
		&\leq& 2\delta+2
\end{eqnarray*}
Thus, if $\delta$ were a-priori bounded, so would be $w(a,b)$ for every special pair. This means --
\begin{cor} Suppose $G$ is a group acting geometrically on a CAT(0) space $X$, and $\HH$ is a $G$-invariant uniform halfspace system on $X$. If $G$ acts co-compactly on $C(\HH)$, then $\HH$ has the strong parallel walls property.\ep
\end{cor}
We would like to find weaker conditions forcing the strong parallel walls property. For example, in the case of Coxeter groups acting on their Davis-Moussong complexes the strong parallel walls property holds for all cases -- not just the co-compact ones. Let us formulate a question:

\vspace{2mm}
\noindent\textbf{Problem: }{\it Find a reformulation of the strong parallel walls property for a (uniform) system $\HH$ in terms of the geometry of the quotient of $C(\HH)$ by $G$. Is it possible that $\HH$ has the strong parallel walls property if and only if the action of $G$ on $C(\HH)$ is geometrically finite in the sense of Wise \cite{[Wise]}? if so, then what can be said about the cusp subgroups?}

\vspace{2mm}
Related to this is the very basic question regarding the dimension of $C(\HH)$:

\vspace{2mm}
\noindent\textbf{Problem: }{\it Is it true that a uniform halfspace system has to be $\omega$-dimensional? finite-dimensional?}

\vspace{2mm}
This question is the reason for the detailed discussion of proper cubings in \ref{subsubsection:proper cubings}. Although it is not true in general that a discrete halfspace system -- even one with a proper dual cubing -- has no infinite transverse subset, it seems possible that uniformness in conjunction with regularity are restrictive enough to at least rule this possibility out.

\bibliography{prop}
\bibliographystyle{alpha}

\end{document}